\documentclass[11pt]{article}
\usepackage{fullpage}
\usepackage{multicol}
\usepackage[usenames,dvipsnames]{color}
\usepackage{hyperref}
\hypersetup{
  colorlinks,
  citecolor=MidnightBlue,
  filecolor=magenta,
  linkcolor=blue,
  urlcolor=red
}

\usepackage[normalem]{ulem}

\usepackage{graphicx}
\usepackage{amsmath}
\usepackage{amssymb}
\usepackage{amsthm}
\newcommand\addtag{\refstepcounter{equation}\tag{\theequation}}

\definecolor{lightgray}{rgb}{.95,.95,.95}

\usepackage{comment}
\usepackage{subfigure}
\usepackage{natbib}

\usepackage{colortbl}

\newcommand{\minimize}{\mbox{\rm minimize }}

\newcommand{\st}{\mbox{\rm subject to }}

\renewcommand{\a}{\mbox{\boldmath $a$}}

\newcommand{\x}{\mbox{\boldmath $x$}}
\newcommand{\s}{\mbox{\boldmath $s$}}
\newcommand{\q}{\mbox{\boldmath $q$}}
\newcommand{\rb}{\mbox{\boldmath $r$}}

\newcommand{\z}{\mbox{\boldmath $z$}}
\newcommand{\y}{\mbox{\boldmath $y$}}
\newcommand{\lambdab}{\mbox{\boldmath $\lambda$}}
\newcommand{\gammab}{\mbox{\boldmath $\gamma$}}

\newcommand{\deltab}{\mbox{\boldmath $\delta$}}
\newcommand{\taub}{\mbox{\boldmath $\tau$}}

\newcommand{\w}{\mbox{\boldmath $w$}}
\renewcommand{\b}{\mbox{\boldmath $b$}}
\renewcommand{\c}{\mbox{\boldmath $c$}}
\renewcommand{\d}{\mbox{\boldmath $d$}}
\newcommand{\p}{\mbox{\boldmath $p$}}

\newcommand{\A}{\mbox{\boldmath $\alpha$}}
\newcommand{\B}{\mbox{\boldmath $\beta$}}

\newcommand{\R}{\mbox{\rm\bf R}}

\newcommand{\Rplus}{\R_+}

\long\def\killtext#1{}
\newlength{\lpskip}
\begin{document}

\begin{center}
	{\sc {\LARGE On representation of energy storage in electricity planning models}}

\medskip
\renewcommand*{\thefootnote}{\fnsymbol{footnote}}
	{\large James H. Merrick\footnote{\emph{Geal Research, County Offaly,
            Ireland.} \texttt{jmerrick@alumni.stanford.edu},
            \texttt{james@geal.ie}}, John E.T. Bistline\footnote{\emph{Electric
            Power Research Institute, Palo Alto, CA 94304, USA.}}, Geoffrey J. Blanford{$^\dagger$}}
\renewcommand*{\thefootnote}{\arabic{footnote}}
\setcounter{footnote}{0}

	\medskip
	\emph{May 31, 2021}

\end{center}

\begin{abstract}

This paper considers the representation of energy storage in
electricity sector capacity planning models. The incorporation of
storage in long-term systems models of this type is increasingly
relevant as the cost of storage technologies, particularly batteries,
and of complementary variable renewable technologies, decline. To
value storage technologies appropriately, a representation of linkages
between time periods is required, breaking classical temporal
aggregation strategies that greatly improve computation time. We
appraise approaches to address this problem, highlighting a
common underlying structure, conditions for lossless
  aggregation, and challenges of aggregation at relevant geographical
scales. We then investigate solutions to the modeling
  problem including
  a decomposition scheme to avoid temporal aggregation at a parallelizable
  computational cost. These examples frame aspects of the problem ripe for contributions from the operational research community.
\end{abstract}


\vspace{20mm}
\fbox{
  \begin{minipage}{0.9\linewidth}
\noindent \textsc{\large \textsc{Geal Research Working Paper Series}}\\

GR-02-v2.
The Geal Research Working Paper Series enables discussion and dissemination of ideas while work is still ongoing. Comments very welcome.

\noindent The views expressed in this paper are those of the authors alone and do not necessarily reflect those of their institutions.

\noindent\textcopyright 2020-2021. This manuscript version is made available under the CC-BY-NC-ND 4.0 license\\
\url{http://creativecommons.org/licenses/by-nc-nd/4.0/}

 \end{minipage}
}
\newpage


\section{Introduction}

This paper considers the representation of energy storage technologies in electricity sector planning models. These models are typically formulated as optimization problems to find least-cost portfolios for power sector investments and dispatch and have been used in the power industry and for associated energy sector policy analysis since the 1950s. They have been extensively discussed in the operational research and management science literatures \citep{Masse1957,Murphy2005,Singh2009}.

Model formulations have been updated in recent years to incorporate
technological developments in, and policy supports for, variable renewables
(especially wind and solar power) and energy storage technologies. In
particular, higher temporal and spatial resolutions are needed to
adequately capture variability, a key economic characteristic for
renewables and storage \citep{Cole_2017,Collins_2017}. Strategies to
capture these features often focus on renewables and not energy storage
\citep{Blanford2018,Merrick2016e}, which entails related but distinct
modeling considerations \citep{Bistline_2020b}. Specifically, temporal
aggregation strategies do not typically include chronology, which is
necessary to represent state-of-charge (energy balance) constraints for
energy storage systems, raising issues of potential interest to the
operational research community.
{Other decisions that depend on chronology
  include unit commitment decisions for individual power
  plants (e.g., considering how many hours a plant might be used) and
  time-shifting demand like electric vehicle charging.}
Reduced complexity of electricity models has large impacts on feasibility,
cost, and emissions outcomes \citep{Bistline_2020c,Bistline_2017b}, which
makes it important to think critically about formulations
(for example, \cite{Bistline_2020} demonstrate how electricity model formulations can alter the sign of emissions impacts from energy storage deployment).

This paper evaluates approaches to address this problem of temporal aggregation in electric sector models with energy storage. Storage technologies have become increasingly important in modeling decarbonization and high renewables scenarios, especially as costs decline and deployments increase \citep{Gorman_2020}. However, storage technologies have complex and diverse cost, value, and performance characteristics that make them challenging to model \citep{Bistline_2020b}. One approach is to use merchant price-taker models with historical market data for hypothetical energy storage systems, but these frameworks omit power system feedbacks and have a limited ability to estimate market depth \citep{Evans_2019,Braff_2016}. Production cost models can assess short-run system operations using detailed simulations of unit commitment and dispatch typically over a year, but only capture static systems without changes in investment and thus cannot estimate market depth \citep{EPRI_2020}. Capacity planning and dispatch models -- the focus of this paper -- can help to assess the long-run value of energy storage by accounting for both investment and dispatch effects over a multi-decadal horizon \citep{Santen_2017}. Temporal aggregation is the focus of this work due to both its importance in driving model results and its operational research dimension; however, it is important to keep in mind that it is only one challenge among many in modeling storage \citep{Bistline_2020b}.
Additional modeling challenges associated with energy storage are representations of technologies (including hybrid systems and cross-sector interactions), market participation, policies/incentives, and spatial aggregation.

The novel contributions of this paper comprise a
conceptual framing of the modeling problem in Section \ref{sec:model},
an illustration of the numerical challenge at the scale of 
large interconnected power systems as exist in North America,
Europe, and China, in Section \ref{sec:thechallenge}, and an
investigation of solutions in Section \ref{sec:solutions}.
A key insight from Section \ref{sec:model} is that aggregation methods
in the literature can be represented as a general representation and
that conditions exist for a general lossless
aggregation. Finding the minimal such aggregation for a given model
input dataset remains
an open problem however. Furthermore, the numerical analysis of Section
\ref{sec:thechallenge} indicates the limited scope for substantial
aggregation in practice across all questions that may be asked of a model. The solutions investigated in
Section \ref{sec:solutions} range from alternative modeling paradigms
to a decomposition
approach that
enables the
solution of large problems,
avoiding the need for aggregation.

\section{Background}
\label{sec:background}

This work investigates the representation of energy storage technologies in capacity planning models, which consider system-level interactions for investment decisions (including storage, generation, and transmission assets) and operational dynamics (which influence and are influenced by investment decisions). Prior work, such as that of \cite{Cruise2019,Zhou2016}, applies operational research approaches to the efficient operation of an individual storage unit.
The models underlying our paper consider different questions, for example,
how much storage is economical to deploy in a given system, how does
storage contribute to public decarbonization policies, and how does it
interact with other technologies across interconnected power systems. More
broadly,
many interconnected decisions, like choosing the level of storage deployment,
are treated endogenously in such capacity planning models, necessitating careful formulation to enable efficient computations while capturing salient system features.

We particularly explore in this paper the
aggregation of model variables and constraints that are defined over
temporal operating periods, and how this aggregation affects the
representation of storage technologies. An example of such a variable is energy discharged from a
battery at each hour.
Temporal aggregation has significant
computational benefits and, in the absence of storage technologies,
was a good candidate for aggregation due to the high amount of
redundancy in associated temporal data. In addition to the
computational benefits, \cite{Merrick2019} discuss how aggregation and
reduced-form representations are conceptually desirable once the
aggregation does not materially distort relevant model outputs.
In the absence of storage, as
\cite{Merrick2016e} discusses, and in accordance with the analysis of
\cite{Rogers1991} and the aggregation bounds of
\cite{Zipkin1980a,Zipkin1980}, a good strategy for temporal
aggregation is the gathering together of similar hours. With the
presence of energy storage however, aggregation methods must also
maintain a representation of the chronology between periods.

\subsection{Value of energy storage}
When we aggregate temporal representation, what representation of energy storage do we wish to maintain? In addition to not distorting the valuation of other technologies, we wish to appropriately value energy storage options, the fundamentals of which we now discuss.

The potential value of energy storage systems is more complex than other technologies due to many services that it can provide and difficulty of capturing all streams in a single framework, since they operate over wide spatial and temporal scales and exhibit location-specific variation based on the grid mix, benefiting parties, and market rules \citep{Balducci_2018}. We focus on value streams that are most relevant in long-run planning models, reflecting projections for potential market depth and services already captured for other technologies; namely, energy and capacity value \citep{Bistline_2020b}.
Representation of energy storage value streams not considered in this paper, such as
ancillary service provision and transmission deferral,
will benefit from the same advances in modeling
  chronology that
  enable representation of the 
  energy and capacity value of storage.

We define energy value as the ability to take advantage of daily, weekly, and
even seasonal arbitrage opportunities, and different technologies are
more and less suited to capture different arbitrage opportunities
\citep{Mongrid_2019}. We define capacity value
as the ability to provide electricity during scarcity
events, when price of power is higher than marginal cost of
production.\footnote{Energy and capacity value are defined
  mathematically in Section \ref{sec:evcv}.} As we next discuss, many
temporal aggregation strategies in the literature make \emph{a priori} assumptions about which value streams matter.

\subsection{Aggregation approaches in the literature}
\label{sec:survey}

This section introduces and discusses numerous aggregation approaches
found in the literature, but does not provide an exhaustive list.
For instance, we are not including the so-called ``infinite reservoir approach,'' which omits the storage balance/state-of-charge constraint. We also do not discuss offline methods, where storage-related assessments use pre- or post-process calculations that can iterate with the main optimization model \citep{Cole_2017}.

\subsubsection{Representative sequences}

Representative sequence approaches, and representative day methods in
particular, are the most common temporal aggregation strategies in the
literature. Representative day approaches harness an intuitive,
interpretable structure in the data: daily cycles. Note that while the
daily cycle structure is readily apparent in load and solar data, in many
regions it is not necessarily present in wind data. The approach no longer
involves the aggregation of hours with similar characteristics, but the
aggregation of days (or other sequences) with similar characteristics,
allowing chronology to be maintained within days, but not necessarily
across them (representative weeks are alternately chosen on
  occasion. For example, see \cite{DeSisternes2013}). The model may choose
to deploy storage units to operate within these representative
days. Papers on how to choose representative days include
\cite{Johnston_2019,Teichgraeber2019,Garcia-Cerezo2019,Liu2017,Nahmmacher2016,Poncelet2016a}. These
papers are largely based on clustering methods, with various adaptations
adopted to improve the choice of representative days. Implementations
typically select 4 to the order of 10 days with hourly resolution or multi-hour blocks,
a sample of which are shown in Table \ref{Table_rep}.

\begin{table}
	\caption{Number of representative days for selected models.}
	\vskip 6pt
	\centering
	\begin{tabular}{l c c}
		\hline\hline
		\multicolumn{1}{c}{\textbf{Paper}}  & \multicolumn{1}{c}{\textbf{Model}} & \multicolumn{1}{c}{\textbf{Rep. Days}} \\
		\hline
		\cite{Jayadev2020} & OSeMOSYS	& 4 \\
		\cite{Nahmmacher2016} & LIMES-EU & 6 \\
		\cite{Despres2017} & POLES/EUCAD & 12 \\                
		\cite{Nelson2012} & SWITCH 1.0	& 24 \\
		\hline
	\end{tabular}
	\label{Table_rep}
\end{table}

These aggregation strategies often do not allow linkages across
representative sequences (e.g., specifying that energy storage
balances must return to zero by the end of each sequence). This focus
on intraday linkages suggests that representative day methods are more
suitable for some energy storage technologies such as lithium-ion
batteries than others. The key assumption however is that days can be
chosen that are representative of the distribution of possible
days. As shown in Section \ref{sec:model}, extreme peak pricing periods matter for energy storage valuation, and clustering methods that choose days that capture the center of the distribution may miss crucial extreme days. In turn, knowing \emph{a priori} which day the peak pricing may fall is challenging, as it is endogenous to the model outcomes and may shift based on the grid mix across different scenarios and regions. \cite{Merrick2016e} shows for a sample dataset, when wind and solar profiles are included in addition to load profiles, the number of unique days increases dramatically. 
In Section \ref{sec:hard}, we illustrate this
  challenge further at the scale of the contiguous United States and
  its associated spatial diversity in temporal profiles.

\subsubsection{System states}

System states approaches identify unique states
  that comprise time periods with similar characteristics and estimate
  a probability transition matrix between states
  \citep{Wogrin2016}. The transition matrix enables the representation
  of chronology.
  To keep track of the hourly energy balance and ensure it is within energy storage capacity, the approach includes a structure to infer the hourly balance from the state transitions. Computationally expensive aspects of this structure are avoided by introducing the assumption \emph{a priori} which hours are `charge hours' and which hours are `discharge hours'. \cite{Tejada-Arango2018} adapt \cite{Wogrin2016}, by combining system states and representative days through modeling short term storage dynamics within the representative period (day), and applying the system state approach to model storage across these periods. In their empirical results, they show their approach performs better than the system states approach alone.

  \cite{Kotzur2018a} present a similar approach
 where representative days are
 designed to represent short term storage dynamics,
while a superposition concept represents seasonal storage dynamics.
The assumption is made that a storage
unit exhibits identical behaviour within each day, just at different absolute
levels. For example, in a world with two days, instead of modeling 48 hours, 24 hours are modeled to represent the relative pattern, and 2 additional variables model the absolute level in each day.

Empirical work indicates (for example calculations in
  \cite{Merrick2016e}) that there are simply a large number
of unique system states to represent, blunting the effectiveness of
the system states
aggregation. However, as discussed in
the context of a more general representation
introduced in
Section \ref{sec:general}, the method theoretically
  can lead to a lossless aggregation.

\subsubsection{Adapted aggregation methods}

\cite{Pineda2018} consider clustering with the constraint that only
adjacent periods may be clustered, allowing the aggregated periods to then
be treated as chronological within a model. This method shows improved
performance relative to choosing representative days and weeks, for the
same number of periods. However, this improved performance is achieved at
a scale in terms of number of periods far greater than most models in the
field, challenging their computational limits (for example,
  \cite{Pineda2018} compare the performance of adjacent clustering
  relative to other methods for 672 aggregated periods (28
  days/4 weeks), a greater resolution than many
    models). In Appendix \ref{sec:adj},
    we see the challenge of the approach at the scale of the contiguous
    United States, with diverse load and renewables profiles across regions.
    If this scale is manageable, it is a promising method to achieve some level of aggregation.

\cite{Zhang2018a} and \cite{Duan2019} discuss state compression of Markov processes, which appears to be a promising framework for reasoning about {general aggregation strategies} that maintain chronology. At the core of the method is singular value decomposition upon the transition matrix between states. 
How many states are required, or how similar members
  of states are required to be, for the representation of storage,
  remain an input to these methods, and depend on the structure of the
  optimization model where the aggregated data are applied.

\section{Framing the Modeling Problem}
\label{sec:model}

We next introduce the core model structure that
  underlies this paper and allows us to frame the
  modeling problem of representing storage, and more generally,
  chronology.
  To
keep the notation compact, the greenfield long-run
minimum-cost model
does not include model variables and constraints tracking
intertemporal dynamics that reflect investments and retirements in
capacity over time.
We also do not include here the spatial dimension representing multiple model
regions and transmission expansion decisions as would be present in many interconnected power systems
models.
Furthermore, the
mathematical model
  structure is this section implicitly assumes perfectly competitive
  markets with foresight and has no representation of stochasticity, a
  feature appropriate for certain questions (for
  example, see \cite{Murphy2005}).

These features and dynamics are important to represent in
many applications and have the sizable secondary effect of increasing
model computational requirements.
That we do not include these features and dynamics in
this mathematical representation does not preclude the associated
discussion and insights from models with these features, as the
  structure represented here is an essence of all such models.
For example, the model we use for numerical calculations later in this
paper, the US-REGEN model \citep{REGEN_2020}, employs the core model
structure outlined here in addition to other features and dynamics.

\subsection{Notation}

\noindent Sets:
\begin{itemize}
  \setlength{\itemsep}{1pt}
  \setlength{\parskip}{0pt}
  \setlength{\parsep}{0pt}
\item $g=1,..,m$ generator types
\item $h=1,..,n$ time periods (e.g. hours)
\end{itemize}
Parameters:
\begin{itemize}
  \setlength{\itemsep}{1pt}
  \setlength{\parskip}{0pt}
  \setlength{\parsep}{0pt}
\item $c^x_g\in\R^n$, variable cost of generation from generator type $g$
  {(typically a vector of constants})
\item $\c^z\in\R^m$, vector of generator capacity costs
\item $c^t\in\R$, cost of storage power capacity (``door'' cost in \$/kW)
\item $c^u\in\R$, cost of storage energy capacity (``room'' cost in \$/kWh)
\item $\d\in\R^n$, vector of electricity demands across time
\item $\a_g\in\R^n$, vector of generator availability
\end{itemize}
Variables:
\begin{itemize}
  \setlength{\itemsep}{1pt}
  \setlength{\parskip}{0pt}
  \setlength{\parsep}{0pt}
\item $\x_g\in\Rplus^n$, vector of generation by generator $g$ 
\item $\z\in\Rplus^m$, vector of generator capacity investment
\item $t\in\Rplus$, storage power capacity (``door'') investment
\item $u\in\Rplus$, storage energy capacity (``room'') investment
\item $\rb\in\R^n$, vector of net charge amount
\item $\s\in\Rplus^n$, vector of storage levels/balance/state-of-charge
\end{itemize}

\subsection{Model}

\[
\begin{array}{rclll}
\minimize_{\x,\z,t,u,\rb,\s} & \sum_{g}c_g^{x}\x_g+\c^{z}\z+c^{t}t+c^{u}u&&&\\
\st &\sum_g \x_g+\rb&= \d&:\lambdab&\\ 
&\x_g &\leq \a_{g}z_g&:\gammab_g&\forall g=1,..,m\\
&s_h&=s_{h-1}+r_h&:\Omega_h&\forall h=1,..,n\\
&|\rb|&\leq t&:\deltab&\\
&\s&\leq u&:\taub&\\
\end{array}
\addtag
\label{mdl:core}
\]

The objective function minimizes the cost of supplying electricity, including the provision of energy storage technologies.
The first constraint, with dual variable $\lambdab\in\R^n$, requires supply in each dispatch period to meet demand. The second constraint, with dual variable $\gammab\in\R^{m*n}$, states that generation of each technology cannot exceed available capacity in each period. The third constraint, with dual variable $\Omega\in\R^n$, tracks the storage balance as the storage unit charges and discharges. Since the current state of the energy storage system is impacted by all previous states (which is uncommon for power sector resources), this state-of-charge constraint requires chronology (i.e., linkages across time) to be represented, which makes storage computationally challenging. The fourth and fifth constraints, with dual variable vectors $\deltab\in\R^n$ and $\taub\in\R^n$ respectively, ensure that a) charge/discharge does not exceed the charge/discharge capacity of the unit (i.e., the power capacity or ``door''), and b) the stored energy does not exceed the storage energy capacity (``room''). Finally, note a number of our variables are defined in the domain of the positive real numbers.

Note that this general formulation can encompass a range of energy storage technologies such as batteries, pumped hydro, and compressed air energy storage \citep{Mongrid_2019}. The cost structure of energy storage is taken as an input, including the power capacity cost ($c^t$ in \$/kW) and energy capacity cost ($c^u$ in \$/kWh). Outputs include the energy storage power capacity ($t$ in kW), which governs the maximum rated charge and discharge rates, and energy capacity ($u$ in kWh), which bounds the total electricity stored. The ratio of energy to power capacity determines the duration that the storage device can provide rated power (or length of time needed to charge) and also can be specified exogenously.

\subsection{Derivations}
\label{sec:value}

\subsubsection{Marginal value from optimality conditions}

For a similar model structure, \cite{Lamont2013} derives optimality
conditions relating to the marginal value of investment in power
charge/discharge capacity (``door,'' kW) and energy capacity
(``room,'' kWh). As noted by \cite{Lamont2013}, the marginal value of
these capacities when they are deployed is the sum of rents on their associated constraints:
\[
c^u=\sum_h\tau_h
\label{eq:room}
\addtag
\]
\[
c^t=\sum_h\delta_h
\label{eq:door}
\addtag
\]

These conditions show that the model will invest in storage capacities until such point that the sum of rents, the marginal value, equals the marginal cost of capacity deployment. All of these values are a function of the grid mix, including the level of energy storage deployment. In an intertemporal model setting, the model will invest such that the net present value of future rents over the time horizon equals the upfront cost of installation.

We next derive some further identities from the model
  structure with the goal of illustrating characteristics of energy
  storage representation we wish an aggregation strategy to capture.
Exploring (\ref{eq:room}) further, and for cases when there is some
positive investment in energy storage capacity, we can
derive the following (See Appendix \ref{sec:proofs} for
  derivations associated with this section):
\[
c^u=\sum_h(\Omega_{h+1}-\Omega_h)^+
\addtag
\label{eq:mvroom}
\]

That is, the marginal value is the sum of the positive differences in the
dual variable of the storage balance constraint.
  \cite{Lamont2013} points out the relevance of cycles in the structure of
  the optimal solution,
where a
cycle comprises the set of periods between the storage unit having zero
energy stored.
  Noting that the point where a monotonic increase in $\Omega$
will reverse is when the quantity of stored energy is at zero, we can derive
the following:

\[
c^u=\sum_k(\max_{h\in k}\Omega_{h}-\min_{h\in k}\Omega_h)=\sum_k(\Omega_{b(k)}-\Omega_{a(k)})
\addtag
\label{eq:mvroom2}
\]

Where the index $k$ is across the set of charge/discharge cycles of a
storage unit, and $b(k)$ is the last period of the discharge portion of a
cycle $k$, and $a$
is the first period of the charge portion of a cycle $k$. While the
location and duration of cycles are endogenous to the
model, the cycle is a useful device for considering what comprises the marginal
value of a storage unit.
Furthermore, (\ref{eq:mvroom2}) has the interesting implication that we can collapse
the marginal value obtained from each charge/discharge cycle into two boundary
prices. 

The associated buying and selling price, $\Omega$,
can be derived as follows,
when the unit is charging or discharging, i.e. $|r|>0$
:

\[
  \Omega_h = \lambda_h+\delta_h^c-\delta_h^d
  \addtag
  \]

  Recalling that $\delta$ is the dual variable (rent) of the door capacity
  constraint, we can see that, when that constraint is not binding,
  $\Omega$ equals $\lambda$, the price of electricity.
$\Omega$ may be interpreted as the local price of electricity facing the
    charging/discharging storage unit, differing from the system price by
    a congestion charge entering / exiting the unit. When buying energy,
    the congestion charge increases the price the unit pays, while when
    selling, the congestion charge decrease the price received.
  The extent to which a unit of storage can claim the full price of
  electricity when arbitraging across periods is thus limited by the price
  associated with charge/discharge capacity. If the charge/discharge
  capacity of a storage unit were free, the marginal
  value of storage would be purely based on electricity price arbitrage.
  Since electricity arbitrage value is dependent on peak and off-peak price differentials, models with limited temporal resolution could dampen price variability and lower the value of energy storage \citep{Diaz_2019}.

As (\ref{eq:mvroom2}) shows, the greater the dispersion in prices, the
greater the marginal value of a storage unit (wind and solar deployment have been shown in modeling studies and market data to increase price variability \citep{Mills_2020}, \emph{ceteris paribus} increasing the marginal value of energy storage). Also, as storage costs decline, it allows the marginal value of the technology to decline while still remaining a viable investment. Note that the technology itself decreases its own marginal value as it is increasingly deployed, linking periods and reducing disparities in prices \citep{Denholm_2019,Bistline_2017,deSisternes_2016,Blanford_2015}. Decreasing returns (``value deflation'') also occurs for wind and solar, as their economic value declines as their penetration increases. These declines have been observed in a range of actual market settings and prospective modeling studies \citep{Wiser_2017,Bistline_2017,Gowrisankaran_2016,Hirth_2013}.

These identities are displayed to provide intuition into the drivers of storage value, and the relationship between energy and power charge/discharge capacity, all at the margin. Such drivers aggregate model representations of energy storage will want to maintain. To analyze these countervailing effects and which dominate (under different conditions), we need a numerical model like the one applied in Section \ref{sec:general} and described in Appendix \ref{sec:regen}.
  At this point, while we have not explicitly modeled features
  like round trip storage efficiency, we note both the role of cycles and the role of dispersion in electricity
  prices are fundamental in the marginal value of a storage unit, indicating these are concepts
  that may be harnessed by aggregation strategies.

\subsubsection{Energy and capacity value}
\label{sec:evcv}
We next consider total value realized by a storage unit deployed in the
model. Energy value and capacity value are a useful distinction for
considering this realized value, and are relevant for assessing the
strengths and weaknesses of various aggregation methods.
Defining $\gammab^*$ as the difference between the electricity price,
$\lambdab$, and the maximum marginal cost for any dispatched technology
(noting from the optimality conditions that
$\lambdab=c_g^x+\gammab_g,\forall g\quad\implies$
$\gamma^*_h=min_{(x_{h,g}>0,g)}(\lambda_h-c_g^x)$), 
we define energy and capacity value as follows:
\begin{itemize}
\item Energy value: Profits in energy arbitrage from a storage unit from prices excluding the scarcity premium in peak pricing periods. We define this as $\rb\lambdab^*$, where $\lambdab^*$ is the electricity price, $\lambdab$, adjusted to remove the scarcity premium, $\gammab^*$.
\item Capacity value: Value of dispatch in scarce periods, i.e. the value
  of the scarcity premium captured by the storage unit,
  $\rb\gammab^*$. Note that the model formulation (\ref{mdl:core})
  omits a reserve margin constraint for concision, so the
  market-clearing constraint with dual variable $\lambdab\in\R^n$
  embodies both the energy and capacity prices, while
   in a model with an explicit reserve margin constraint, the dual of
   that constraint will carry the capacity price.
\end{itemize}  

We would like any representation of energy storage in a model to allow
both these value streams to be realized, as missing either can lead to
misvaluation. For these streams to be realized, electricity prices,
including scarcity premiums, need to be represented and adjusted
appropriately as the model endogenously chooses storage deployment levels
and the capacity mix of other technologies. Additionally, restricting
storage deployment to fixed ratios of solar deployment, for example,
considering only solar plus storage as a joint technology, could miss some
of the associated value of either individual technology as prices
endogenously adjust to different model outcomes.
As discussed in Section \ref{sec:hard}, capacity value in particular requires a representation of peak pricing periods, which are difficult to identify \emph{a priori} for large interconnected multi-region systems with many generation technology options.

\subsection{General representation of an aggregated model}
\label{sec:general}

This subsection develops a generalized formulation of a capacity planning
model with energy storage that encapsulates both the non-aggregated
formulation and aggregated approaches discussed in Section
\ref{sec:survey}.
This formulation illustrates common features, and common
strengths and weaknesses, across aggregation methods, with a view to
aiding design of future improvements.
With the introduction of the following parameters, we can develop a more general version of (\ref{mdl:core}) from Section \ref{sec:model}: 
\begin{itemize}
\item $\w\in\R^n$ weight associated with each period
\item $P\in\R^{n\times n}$ state transition matrix, where $P_{ij}$ is the probability of jumping from any period $i$ to any other period $j$
\item $\q\in\R^n$ duration in a state before transitioning to another state
\end{itemize}

Introducing $\w$ into the objective function, and $P$ and $q$ into the storage balance constraint, allows us to present our problem's model structure as in (\ref{mdl:agg}) below:
\[
\begin{array}{rcll}
\minimize_{\x_g,\z,t,u,\rb,\s} & \w\sum_{g}c_g^{x}\x_g+\c^{z}\z+c^{t}t+c^{u}u&&\\
\st &\sum_g \x_g+\rb&= \d&\\ 
&\x_g &\leq \a_{g}z_g&\forall g=1,..,m\\
&\s&=P\s+\q\rb&\\
&|\rb|&\leq t&\\
&\s&\leq u&\\
\end{array}
\addtag
\label{mdl:agg}
\]
This formulation is equivalent to our earlier formulation (\ref{mdl:core}) when (i) $\w$ consists of the vector of ones, (ii) $\q$ consists of the vector of ones, and (iii) $P$ is the identity matrix with all ones shifted one place to the right, i.e., $P_{h-1,h}=1$ for all $h$.

Furthermore, the aggregation methods mentioned above, with a compressed temporal dimension, may all be considered as special cases of the formulation (\ref{mdl:agg}). For example, in a representative day scenario, $\w$ encodes the weight of each day, $\q=1$, and $P$ encodes the representative day structure with cycles of 24 states each, connected or unconnected depending on how intraday storage is treated. For the adjacent hierarchical clustering proposal of \cite{Pineda2018}, $\w=\q$ and $P$ comprises again the identity matrix shifted one unit to the right.

\begin{figure}
  \begin{minipage}{0.45\textwidth}
  \centering
  \includegraphics[width=0.97\textwidth]{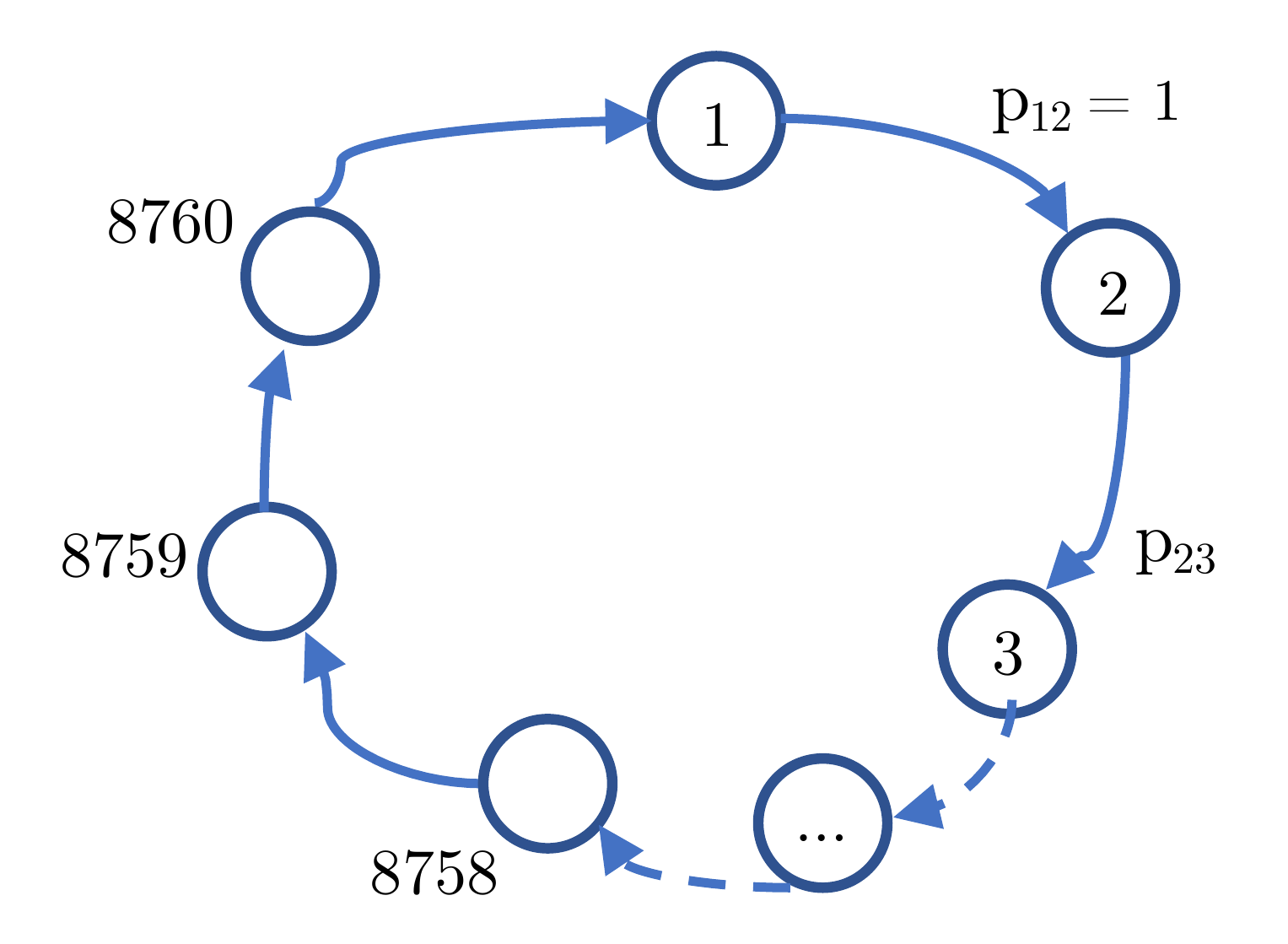}
\end{minipage}\hfill
\begin{minipage}{0.45\textwidth}
  \centering
  \includegraphics[width=0.97\textwidth]{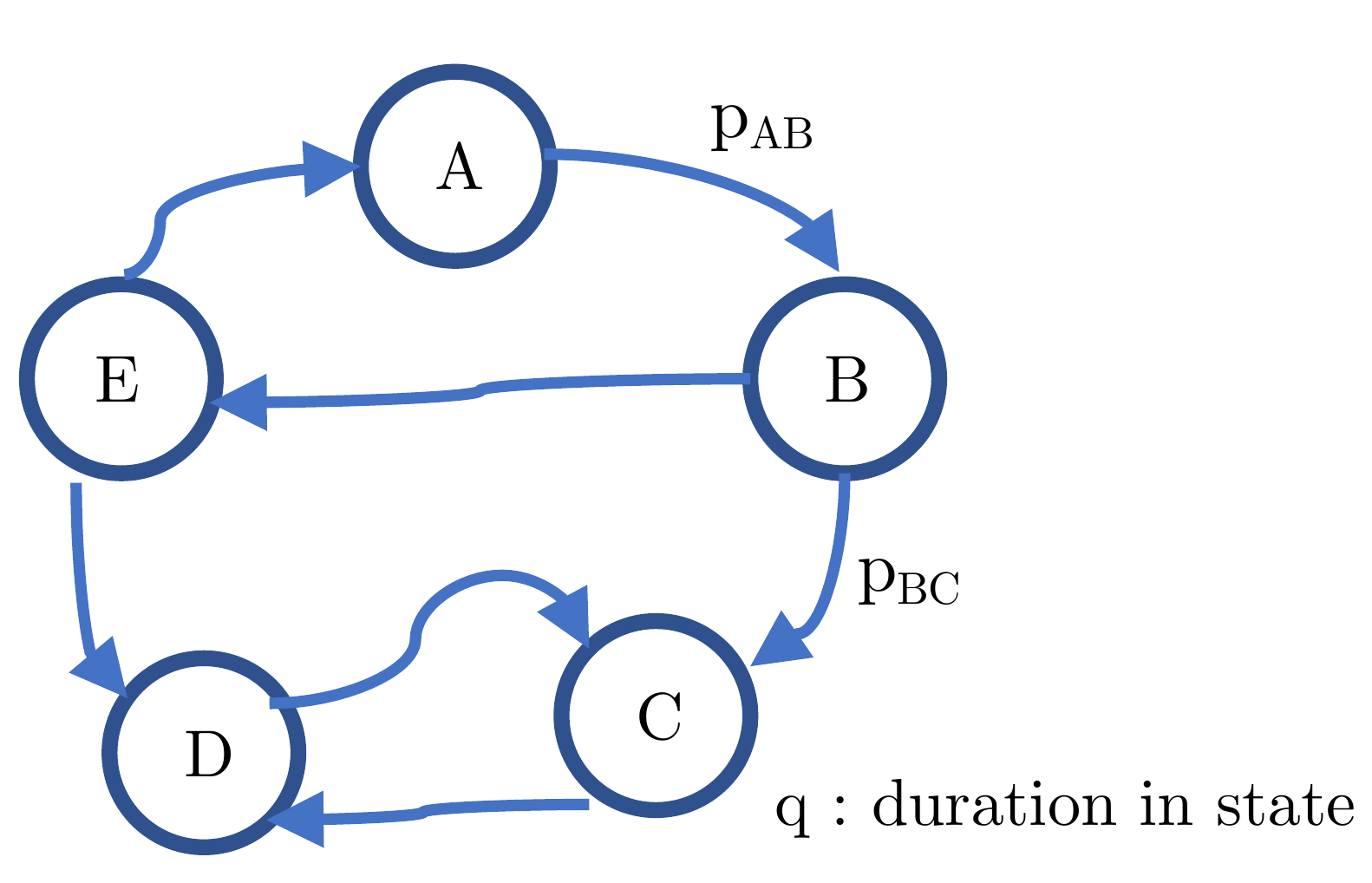}
\end{minipage}
\caption{Aggregation of hourly data into system states. For example, let us say hours 1, 25, and 39 are similar to each other and then could be aggregated to a common state, State $A$.}
\label{fig:state_transition}
\end{figure}

System states in its raw form involves aggregating hours and then, given the aggregation, calculating $P, \q, \w$ empirically for that aggregation. $P, \q, \w$ provide a synthetic mapping to represent reduced-form chronology. Figure \ref{fig:state_transition} illustrates how each of the hours can be grouped together with other hours that are similar into aggregate States $A,B,C,\dots$ with the probabilities calculated empirically of moving from one state to another or remaining in a particular state. Let us say that the resulting $P$ matrix implies the stored energy incoming to State $C$ is 0.25 from State $B$, and 0.75 from State $D$. The associated interpretation is challenging, with a plausible interpretation that each state has an ``expected value'' of storage available when making dispatch decisions in the context of a planning model.

Similarly, the optimality condition on the marginal value of energy storage capacity (\ref{eq:mvroom}) shifts, under the more general formulation, as follows:
  \[
\quad c^u=\sum_h(\Omega_{h+1}-\Omega_h)^+ \quad \rightarrow \quad c^u=\sum_{{s}}(\p_{{s}}{\mbox{\boldmath $\Omega$}}-\Omega_{s})^+
  \]
This states that the value the energy storage unit receives for transferring energy to the future shifts from the price of the next period to the \emph{expected price} across all periods to which the current period is linked.
    Mapping this identity to cycles is more abstract in this setting, and
    can aversely affect the performance of
  straightforward implementation of this method, which
    is computationally explored in Section \ref{sec:mv}.

\begin{figure}
  \centering
  \begin{subfigure}
  \centering
  \includegraphics[width=0.5\textwidth]{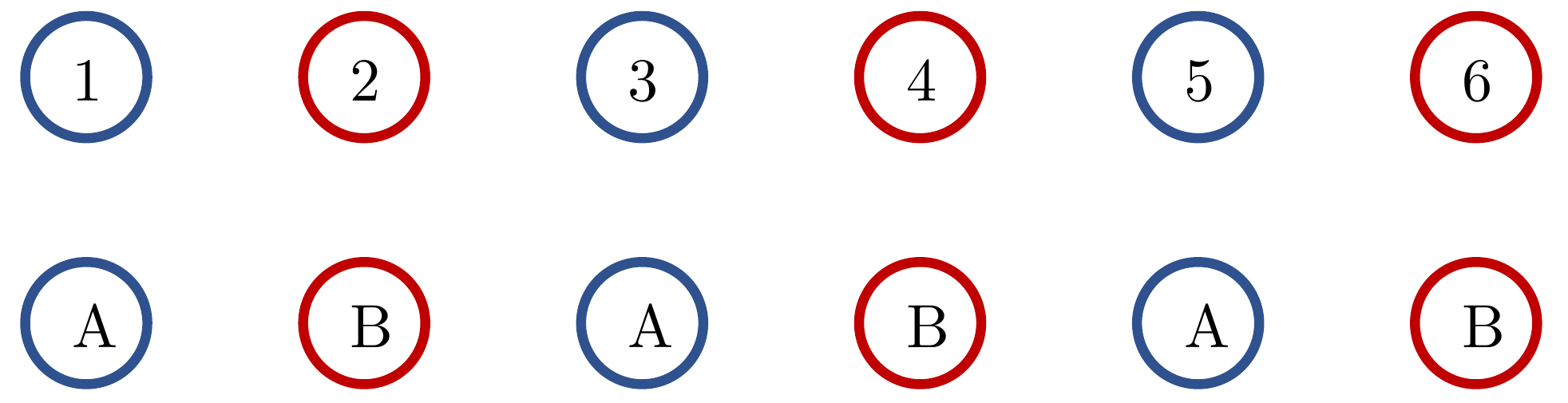}
  \caption{Illustrative example: periods classified into like states}
\label{fig:qw1}
  \end{subfigure}

\begin{subfigure}
  \centering
  \includegraphics[width=0.6\textwidth]{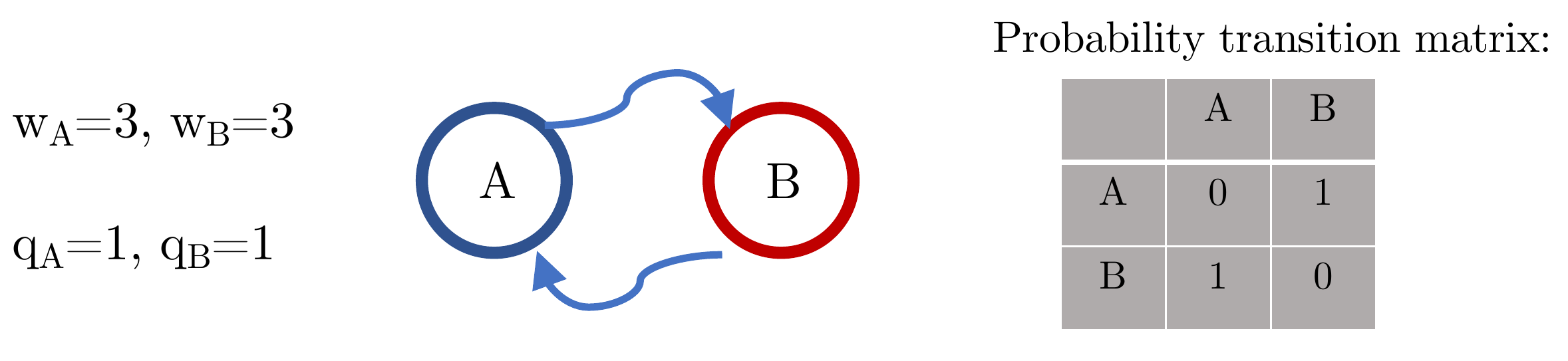}
  \caption{Illustrative example: aggregation of Figure \ref{fig:qw1}}
\label{fig:qw2}
\end{subfigure}
\end{figure}

Reflecting on Figure \ref{fig:state_transition}, it might appear that a `raw' system states approach could never do well, allowing unrealistic energy transfers. Figures \ref{fig:qw1} and \ref{fig:qw2} illustrate a thought experiment where, by inspection, we can see that an aggregated model is equivalent to a non-aggregated model. In this example, the number of states is compressed from six to two while maintaining a representation of chronology. More particularly, in this case, our six states comprised two states repeated in a certain order, an order where it was possible to capture the chronology. Crucial also is the recognition of the distinction between the weight of a state in the objective function, $\w$, and how long one remains in a state upon entering, $\q$. 
In our thought experiment, due to structure identified in
the input data, it is possible to parameterize the general
formulation (\ref{mdl:agg}) in such a way that aggregation does not distort model outcomes. 

We can generalize from this experiment to show (see
  Appendix \ref{sec:proof-ident}) that an aggregation scheme comprising a
mapping between hours $h$ and states $s$ that can retain the following properties will allow an aggregated model be equivalent to a non-aggregated model:

\begin{itemize}
\item Hours that are members of a given state have equivalent temporal characteristics:
      $a_{i}=a_{j}$, $d_{i}=d_{j}$ for all hours $(i,j)$ mapped to a state $s$
\item Each state has only one incoming connection to another state, and
  only one outgoing connection: each row of the $P$ matrix contains one
  $1$ entry, with the remainder all zeros
\item For each $s$, each subsequence of hours that maps to it is
  of equal length, $q_s$
\end{itemize}

As an example, a year comprising 365 identical days would fit these
conditions. Similarly, a year with two unique days, one always following
the other, would also meet these conditions. Additionally, the previously
mentioned adjacent clustering idea would also meet these conditions. Note
that while this would be an aggregation with a guarantee, it is not proven
as a minimal representation. However, we may ask is there a systematic way
of aggregating model input data in such a way as to meet these conditions,
particularly when structure is not obvious in real data. Also, an
aggregation that is simply close to meeting these conditions might be
sufficient for real world applications. We pose this question for future
research.

\section{The Numerical Challenge}
\label{sec:thechallenge}

  Following our conceptual framing of the modeling problem, this
  section illustrates the numerical challenge to represent energy
  storage technologies at the spatial and
  temporal scale of models
  used for pressing policy and strategic questions.

\subsection{Challenge of identifying critical days}
\label{sec:hard}

As discussed in earlier sections, one feature of the capacity planning problem with energy storage that makes temporal aggregation challenging is that we do not know \emph{a priori} which extreme periods are the important ones in driving energy storage's value. In part, this difficulty arises from the dependence of prices on the grid mix and level of storage deployment. This endogeneity suggests that a temporal aggregation method should work for arbitrary levels of energy storage deployment, variable renewable shares, and other system conditions (i.e., capture all relevant extremes of time-series variables such as hourly load, potential wind output, and potential solar output).

To propose an approach to capture these extremes and to illustrate the challenge of representative day selection, this subsection conducts an empirical investigation with time-series data from the United States that populate the Regional Economy, Greenhouse Gas, and Energy (REGEN) model \citep{REGEN_2020,Bistline_2019,Blanford2014}.
The model comprises an intertemporal planning model for the contiguous
United States that chooses electricity capacity investments and
retirements over time that meets electricity demand in line with
modeled policies and technical constraints. Its core
  logic is consistent with the model structure introduced in Section \ref{sec:model}.

\cite{Merrick2016e} showed that the number of unique days in a sample dataset was in the hundreds, indicating the magnitude of the challenge of finding an appropriate aggregation based around representative days, given the aforementioned difficulty of knowing \emph{ex ante} where the peak pricing periods may lie. 
Referring to the previously introduced concepts of energy value and capacity value, both are endogenous to the dispatch \emph{and} capacity choices made by the model. We do know that the peak price that drives capacity value will occur when the system is under strain, for example when demand is high and renewables availability is low. Even more straining on the system is when there are a sequence of consecutive hours with such properties, since storage systems are energy-limited resources. 
The other extreme of low prices also matters for valuation of variable renewables and storage. These occur for example when demand is low and renewables availability is high.

The idea behind this experiment is to compress an annual dataset of temporal availability across the United States into `cumulative days', summing load and variable renewables within each day, before normalizing. Over the resulting 365 days in dataset we apply the strategy of \cite{Blanford2018} to search for selecting `extreme days'. This procedure involves working out the day closest to each load-wind-solar extreme vertex and choosing the minimum number of days such that, for each region, at least one day is chosen that it is within some radius of each vertex point.

\begin{figure}
  \centering
  \includegraphics[width=0.75\textwidth]{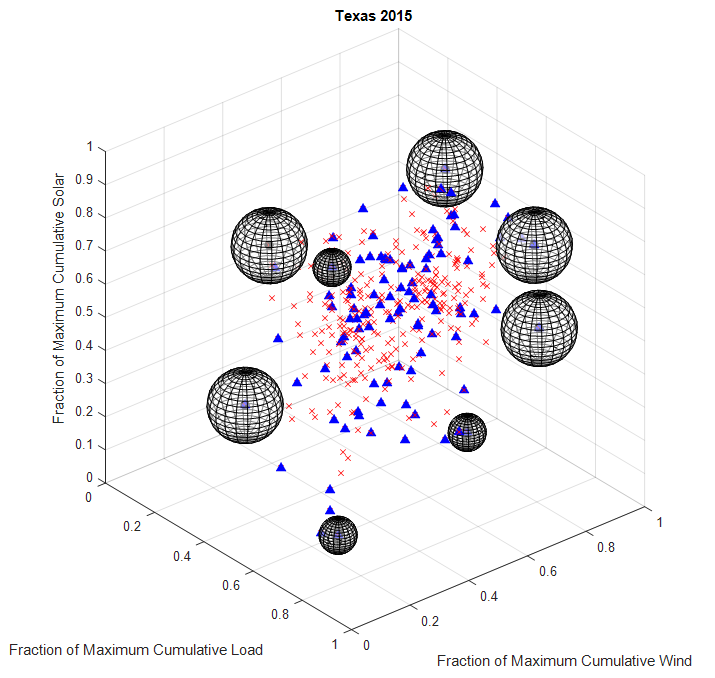}
  \caption{Extreme cumulative days example for Texas using normalized load, wind, and solar data from the US-REGEN model \citep{REGEN_2020}. Red points are cumulative days (365 in the annual dataset), the bubbles are the local vertex requirements, and blue points are the chosen days. Blue points inside a bubble are the days that meet the local requirements, and blue points outside the bubble meet the requirements in other regions.}
  \label{fig:cumulative_day}
\end{figure}

Figure \ref{fig:cumulative_day} displays the results for one region
(Texas) in the 15-region dataset. The extreme hour finding algorithm
in this experiment returns with circa 100 extreme days to meet these
requirements across the contiguous United States, not 4 to 10 days as
is typically used in ``representative hour'' implementations for
electric sector models (Table \ref{Table_rep}). This result highlights the challenge of finding a small set of representative days. Note that this is a \emph{lower bound} on the number of days required, as it ignores diverging patterns within the days, along with the necessary days to cover interior points of the distribution, which could be important in capturing storage value during normal periods.

Also note that this experiment only considers three time-series variables:
load, potential wind output, and potential solar output. If you are in a
region where one renewable resource is known \emph{a priori} as not in-the-money, then the problem may be easier. However, there could be additional time-series variables that may correlate with energy storage value (e.g., hydro output) or additional series to represent wind/solar technologies with different resource profiles, which would increase the number of representative days. It is also important to note that, although the number of required days could be lower if only a single region is considered, the interconnectedness of markets is important for capturing system operations and economics, which increases the dimensionality of the optimization problem. It could also be desirable for a single-region application to have many subregions to characterize differences in resources, which would increase the number of representative days.

\subsection{Performance of selected methods}
\label{sec:mv}

To learn about the in-practice performance of a number of selected
aggregation methods, we conduct a series of numerical experiments with a
`static' version of the aforementioned US-REGEN model. In the static
version, we only solve for one model year, not the full intertemporal
problem. The reduced computational burden allows us to run the
non-aggregated 8760 resolution case, and compare associated model outputs
to model outputs with aggregated resolutions. A challenge to conduct
meaningful experiments in this fashion is that while a given aggregated
method may perform well under the certain set of circumstances of a
numerical run, it may not do so well under different model parameterizations.

Table \ref{tbl:runsA} compares outputs from a) the 8760 non-aggregated
case, b) what we term the `Expected Value Method',
direct implementation of the general formulation outlined in
Section \ref{sec:general}, and c) a representative day formulation. For the
latter, we choose the peak and median load days from each season of the
year (8 days, 192 periods), while for the former we choose 192 periods in
the fashion described by \cite{Blanford2018} that preserves load, wind, and solar
distributions, and directly calculate a probability transition matrix
between periods. The scenario is a carbon price scenario designed to test
the performance of the methods in a case where variable renewables are
incentivized, and the static model year is 2030. To focus on the impact of
renewables and storage, there is a constraint on new
nuclear construction, but this also leads to a narrowing of model options
under the carbon price scenario, which in turn
reduces the scope for error from an aggregation method, highlighting the
challenge of numeric assessment of methods that we wish to perform generally.

  Perhaps the most striking takeaway from Table \ref{tbl:runsA} is the dependence on the question
  being asked of the model when evaluating aggregation methods. We see
  that the aggregation only instigates a small change in the objective
  value of the model, while creating a dramatic speedup in runtime (RAM
  constraints lengthen the 8760 runtime in this example). By many
  criteria, the change in objective function may be considered a very
  small price for the computational gains. However, if the question is what
  is the market for storage capacity under this carbon price
  scenario, we see a greater deviation (the `net capacity value' is
    higher for the Expected Value method in this case as it deploys less capacity,
    keeping prices higher, allowing the capacity it does deploy to achieve
    greater returns).
  For the sort of
  applications these models are used for, and the questions asked, we
  often want whole subcomponents of the aggregated model to match
  corresponding components of the original model. Similarly Table
  \ref{tbl:runsA} contains results spatially aggregated across the USA, masking more
  divergence from the 8760 case at regional level.

  Another feature of Table \ref{tbl:runsA} is that given the nature of
  representative day choice, with many potentially
  critical system days omitted, how does the representative day aggregation do so well? We know that days
  have an intrinsic structure and that if load
  is still the dominant effect, a
  representative day choice based on load, as this one is, could perform well.
Similarly, we consider lithium-ion battery storage, and not longer
cycle storage forms.
  As regards renewables representation, perhaps errors in one region are cancelling out errors in others.
  However, there may also simply be an element of luck in this case.
  Table
  \ref{tbl:runsB} shows a repeat of the experiment with a slight change in
  how the representative days are chosen. Selection B and C respectively
  replace the median load day in each season with the day either side of
  the median in the ranking of load days.
  With this slight change, we can see
  divergence in outcomes, significant in some cases, and this could be
  amplified in different scenarios.

\begin{table}
	\caption{{Comparison of aggregated model outputs relative to
          non-aggregated 8760 case. The scenario is a carbon price scenario
          (\$81/tonne) for the 2030 model year with assumptions about limited CCS availability}}
	\vskip 6pt
	\centering
	\begin{tabular}{l | c | c c}
		\hline\hline
		\multicolumn{1}{c}{\textbf{Model Output}} & \multicolumn{1}{c}{\textbf{8760}} & \multicolumn{1}{c}{\textbf{Expected Value}} & \multicolumn{1}{c}{\textbf{Rep. Days}} \\
		 & & \multicolumn{1}{c}{\textbf{(relative)}} & \multicolumn{1}{c}{\textbf{(relative)}} \\
                \hline
		Storage Room (GWh) & 98.33 & 0.96 & 0.23 \\
                Storage Door (GW) & 37.42 & 0.6 & 0.37\\
		Objective Value (Billion \$) & 176.386 & 1.016 & 1.004 \\          
		Net Capacity Value {(\$/kW)} & 15.65 & 2.62  & 0.98\\
		Net Energy Value {(\$/kW)}& 20.26 & 0.75 & 0.55\\
		CO$_2$ Emissions (MtCO$_2$) & 370 & 1.06 & 1.13\\                                
                Variable renewable energy [VRE] (TWh) & 1627 & 0.99 & 0.98\\
                Natural Gas Capacity (GW) & 434.05 & 1.02 & 0.94 \\
                Speed (s) & 6504 & 0.004 & 0.003\\                
		\hline
	\end{tabular}
	\label{tbl:runsA}
\end{table}

\begin{table}
	\caption{{Comparison of varying representative day choices}}
	\vskip 6pt
	\centering
	\begin{tabular}{l | c c c}
		\hline\hline
		\multicolumn{1}{c}{\textbf{Model Output}} &
                \multicolumn{1}{c}{\textbf{Rep. Days A}} &
                \multicolumn{1}{c}{\textbf{Rep. Days B}} &
                \multicolumn{1}{c}{\textbf{Rep. Days C}} \\
		 & \multicolumn{1}{c}{\textbf{(relative)}} & \multicolumn{1}{c}{\textbf{(relative)}} & \multicolumn{1}{c}{\textbf{(relative)}} \\
                \hline
		Storage Room  & 0.23 & 0.59&0.7\\
                Storage Door  & 0.37&0.81&0.74\\
		Objective Value  & 1.004 &0.94&0.9 \\          
		Net Capacity Value  &  0.98 &1.17&0.91\\
		Net Energy Value &  0.55 &0.39&1.21\\
		CO$_2$ Emissions &  1.13 &0.72&0.8\\  
                VRE  &  0.98 &1.24&1.12\\
                Natural Gas GT Capacity  & 0.99 &0.99&0.99\\
                Natural Gas CC Capacity  & 0.91 &0.74&0.74\\                
		\hline
	\end{tabular}
	\label{tbl:runsB}
\end{table}

\begin{figure}
	\centering
		\includegraphics[width=0.99\textwidth,trim={1cm 4cm 1cm 4cm},clip]{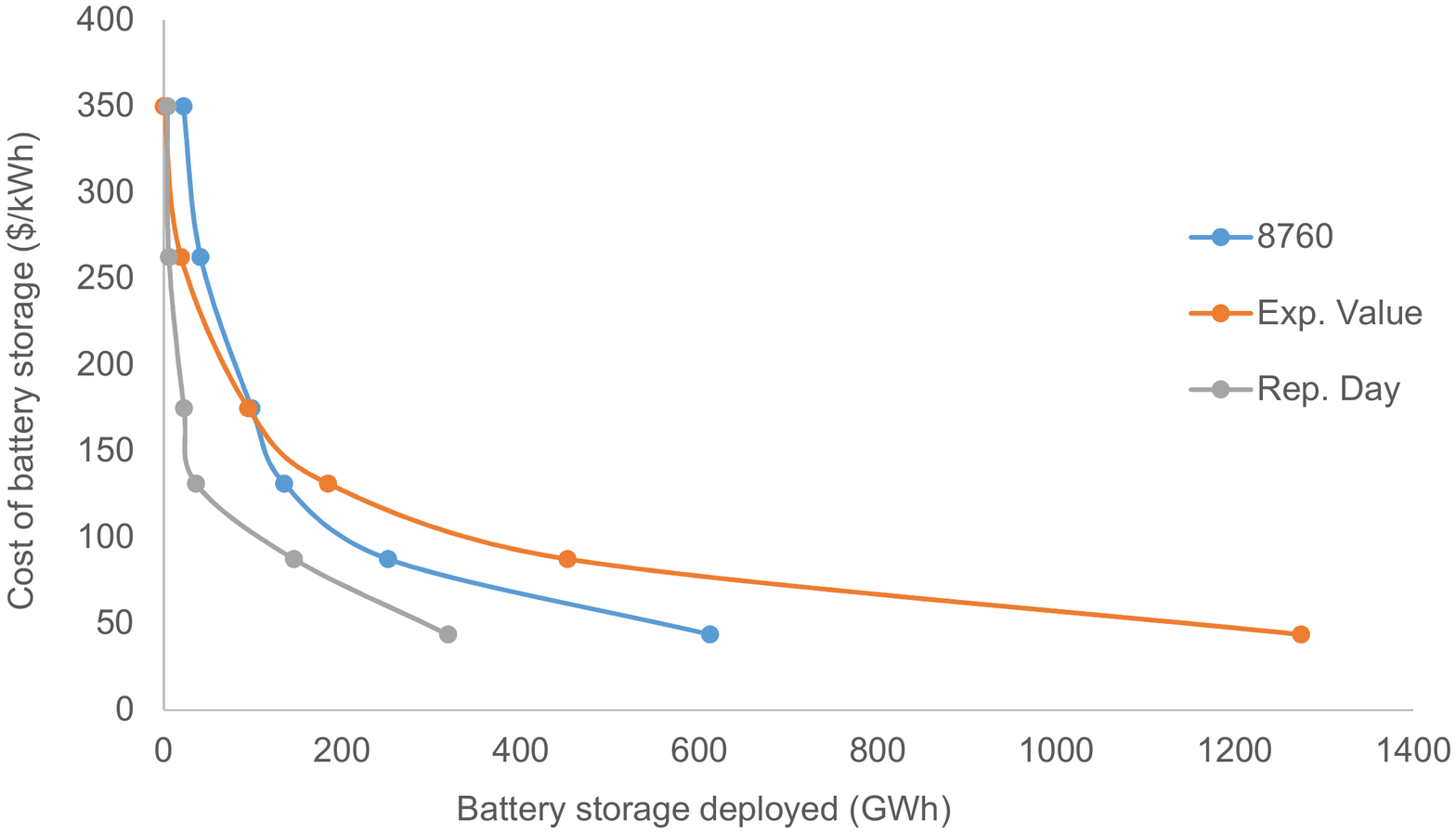}
		\caption{{Marginal value curve of storage across different
                  temporal representations}}
		\label{fig:mv}
\end{figure}

Figure \ref{fig:mv} plots deployment of storage `room' capacity, across
different costs, showing a marginal value curve. Note again that this is for
only one instance of the parameterization. Table \ref{tbl:runsA} referred
to the \$175/kW point on the y-axis.
We see the `expected value curve' matching the `8760 curve' along a portion of
the curve and the `representative day curve' doing better on another.
At the lowest storage cost, we see the representative day method
undervaluing storage, likely missing opportunities due to days omitted,
while the expected value method overvalues storage, seeing unrealistic
opportunities in its probabilistic connections between states.
As discussed, how much
these deviations matter depends upon the question and applications at hand.

\section{Solutions}
\label{sec:solutions}

Based on our exploration of the modeling problem of
  representing energy storage, this
  section outlines solutions, and pathways to solutions, to the
  representation problem, particularly,
  a) improvement of aggregation methods, b) alternative modeling paradigms, and c)
  avoiding aggregation through the harnessing of computational and
  algorithmic advances.

\subsection{Improvement of aggregation methods}

In the improvement of aggregation methods, the goal is to find temporal aggregation strategies that are better than excluding chronology and/or energy storage technologies from long-term models, since omitting energy storage is likely not an acceptable approach given current costs and expectations for future competitiveness. In addition to speed/feasibility and the capacity to address intertemporal questions, other desirable properties in a temporal aggregation strategy relate to the accuracy of the approximation, including the ability to:
\begin{itemize}
	\item Characterize both energy and capacity value accurately
	\item Reflect changes in the marginal value of energy storage at different levels of deployment (namely, decreasing marginal returns at increasing penetration levels)
	\item Retain chronological information across diurnal, weekly, and seasonal timescales
	\item Preserve joint variation of all hourly time-series data across all model regions (e.g., so that variable renewables are simultaneously valued accurately)
\end{itemize} 
Thus, a reduced-form representation should work for arbitrary levels
of energy storage deployment, variable renewable shares, and other
system conditions. As an example, Tables \ref{tbl:runsA} and
  \ref{tbl:runsB} showed an example of aggregation strategies not
  robust to different model parameterizations.

Section \ref{sec:general} showed that the aggregation approaches from the literature
introduced in Section \ref{sec:survey} could all be considered instances
of a general aggregate representation. Furthermore, conditions were
shown where this general aggregation would produce the same model
outputs as a non-aggregated model. Of the methods discussed, adjacent clustering met
these conditions, however further compression was possible through
identification of cyclic structure in model data. A systematic
approach that could find the most compressed representation that meets
the conditions, would improve all the methods discussed, and
essentially collapse them to one.

\subsection{Alternate paradigms}

Underlying our challenge of aggregating electricity planning models in the presence of storage is that we are, for a particular set of assumptions, trying to find the minimum-cost solution across time and all possible capacity mixes that conforms to policy and technical constraints (while evaluating investment and dispatch simultaneously). An alternate to finding aggregation strategies or greatly increasing computation time is to consider two inter-related options of (a) changing the question, or (b) changing the modeling paradigm.

As regards the question, strategies that are employed include a strong decoupling of planning and operation models, where storage is only thoroughly represented in an operation model that considers one year as a snapshop, or alternately, a planning model considers one future year only \citep{Bistline_2020,Brown_2018} with all 8,760 hours and annualized costs. However, the challenge with both these approaches is that, with technological development like energy storage technologies, the planning and operations modes have become increasingly coupled, with outcomes in one strongly affecting the other. Although a ``static analysis'' approach allows many insights and is well-suited for certain questions, it misses intertemporal dynamics that are important for evaluating the time path of investments and retirements, which include meeting intermediate goals with long-lived capital investments.
A related method is to use a ``sequential myopic'' (sometimes called ``recursive dynamic'') approach, where each year is solved individually and the capital stock is carried over. However, this strategy  also misses the intertemporal/forward-looking nature of the capacity planning and dispatch problem.

As regards the modeling paradigm, we can note that we have seen the
capacity value rest on a small subset of periods, which we do not
necessarily know \emph{a priori}. This could perhaps be harnessed by
solving the model multiple times, starting with a coarse resolution,
learning from the solution, and adding in additional periods as
warranted, an approach similar in philosophy to column
generation. \cite{Munoz2016} and \cite{Teichgraeber2020} present
approaches of this type. This broad idea, along with its interactions
with model treatment of uncertainty, is a promising area for future
research.\footnote{For example, does uncertainty in model input data
    warrant simplified model structures?}

\subsection{Avoiding temporal aggregation}
\label{sec:decomp}

Decomposition approaches may allow the harnessing of
  computational and algorithmic advances to avoid temporal aggregation
  altogether. This section outlines a scheme for applying the Alternating Direction Method of Multipliers (ADMM) approach to capacity planning problems at large scale for intertemporally consistent decisions.
\cite{Hoschle2018} propose an ADMM-based method for computing risk-averse
equilibrium in electricity markets, whereas, with the exception of \cite{Frew2016} there has been limited
applications of ADMM to our knowledge for capacity planning problems.

The ADMM approach, discussed further in \cite{Merrickdiss} and
\cite{Merrick2017}, and drawing upon  \cite{Bertsekas2015,Boyd2011,Bertsekas1989,Ye2020}, decomposes the general problem of minimizing $f(\x)$ subject
to $Ax=\b$ into $n$ blocks as follows:

\[
\min f_1(\x_1)+\dots+f_n(\x_n)
\]
such that
\[
A_1\x_1+\dots+A_n\x_n=\b\quad\text{or}\quad A_i\x_i-\y_i=0\text{ }\forall i, \sum_{i=1}^n\y_i=\b
\]

The $\x_i$ blocks are then solved separately, and allowably in parallel:
\[
x_i^{k+1}=\text{argmin}_{\x_i} f_i(\x_i)-\alpha_i^k(A_i\x_i-\y_i^k)+\frac{\beta}{2}||A_i\x_i-\y_i^k||^2
\]

After solving for each block of $\x_i$ variables, the dual variable vector, $\A$ and `target' vector $\y$ are updated, with the magnitude of the update depending on the global $A\x-\b$ residual. The process is then repeated until the residuals converge to $0$. The bulk of the computational work is done in solving the primal $\x$ values at each iteration, and the ability to solve this in parallel allows harnessing of distributed computing to solve a problem at scale. The updating of the dual variables $\A$ and $\y$ are closed-form updates on a central computer.

In the context of our electricity sector planning model application, an important consideration is how to choose the allocation of variables to different blocks. A major consideration is how many constraints are `split' by the decomposition. Each constraint that is split adds to the dimensionality of the convergence criteria, and the number of dual variables and `target' variables to be updated at each iteration. 
However, not splitting enough constraints can imply, depending on the structure, too many blocks for efficient convergence \citep{Zhu_2019}.

There are also useful economic interpretations of the decomposition. We have implemented this scheme in the aforementioned US-REGEN model by  separating operation decisions and capacity decisions, 
with system wide consequences of decisions priced in
  through the augmented terms in the objective function.
Furthermore, operation decisions are solved into separate fast solving hourly dispatch decisions.
Price and quantity information are traded back and forth between the individual dispatch blocks and the investment/retirement capacity decisions at each iteration.

While convergence is guaranteed, as has been noted in other
applications of ADMM, the number of iterations to convergence and
solution quality can be a challenge in practice \citep{Boyd2011}, so
the appropriateness of use of this method depends upon the question
being asked. A further challenge is that manual implementation of the
algorithm may be required (as was the case in US-REGEN implementation), rather than simply handing off a specified model to a solver. It is not yet a catch-all solution, however, but it does provide solutions to problems previously unsolvable, including the ability to evaluate a benchmark to check the accuracy of temporal aggregation methods for intertemporal capacity planning and dispatch models with energy storage.

\section{To Conclude}
\label{sec:conclusion}

This paper has provided context to improve the representation of
energy storage in electricity planning models, and thus indirectly
improve the modeling that informs important societal decisions about
the power sector, the energy sector more broadly, and decarbonization
strategies. The goal of this paper has been to consolidate work on
aggregation for electricity sector planning models with energy
storage, and to build a rigorous foundation for future modeling
progress.
Specifically, we contribute to this literature by
a) appraising
approaches to address temporal aggregation in electricity planning
models,
b) framing the modeling problem of developing a general representation of
an aggregation that values storage technologies appropriately, c)
illustrating from a 
numerical perspective the challenge of finding such an aggregation at relevant
geographic scales, and d) investigating solutions to the problem.
The core challenge is identifying methods that can represent
chronology while being robust to the wide
variety of technology and policy scenarios a planning model of the
electricity sector must consider.
We also note that we frequently want an aggregated model robust in
  representing all model outputs and not only the aggregate
  objective function value.

A substantial increase in periods appears unavoidable within 
the
realm of existing modeling methods, and this work ideally will enable clarity when proceeding so that modeling support for large-scale public and private electricity related decisions are not distorted simply by model structure. In this way, our analysis extends the literature demonstrating how the reduced complexity of electricity models has a large influence on feasibility, cost, and emissions outcomes \citep{Bistline_2020c,Diaz_2019,Blanford2018,Mallapragada_2018,Brown_2018,Bistline_2017b,Cole_2017}.

Section \ref{sec:solutions} identified a number of questions and
directions for future research, including a number
particularly ripe for contributions from the operational research community.
Furthermore, while the focus in this paper is
detailed power systems models, there are research implications for
related models. For example, reduced-form representations of electric sector planning are embedded in integrated assessment models with broader geographical and sectoral coverage, which comes at the expense of limited spatial, temporal, and technological detail \citep{Santen_2017}. Integrated assessment models do not capture storage operations or investment explicitly but assume deployment based on solar/wind penetration. Given how we show here that \emph{ex ante} assumptions about storage resources is likely an overly restrictive framing, how can detailed power system models be used to inform higher-level models that attempt to capture investment and dispatch dynamics? Conversely, how can more detailed energy storage system models and production cost models be used or linked to inform capacity planning models \citep{Bistline_2020b}?


\newpage

\section*{Acknowledgements}
The authors would like to thank those who have read
previous versions of this manuscript, as well as participants in two workshops, for their many helpful suggestions. The views expressed in this paper are those of the authors alone and do not necessarily reflect those of their institutions.

\bibliography{GR-02-v2_arxiv}
\bibliographystyle{apalike}

\newpage
\appendix
\section{Identities from Optimality Conditions}
\label{sec:proofs}
\subsection{Marginal value of storage}
We may write the Lagrangian of the model
  (\ref{mdl:core}), with $\beta$ the vector of dual variables associated with the
  non-negativity constraints,
  as follows:
\[
  \begin{array}{rl}
    \mathcal{L}=&\sum_{g}c_g^{x}\x_g+\c^{z}\z+c^{t}t+c^{u}u\\
    -&\lambdab.\left( \d-\sum_g\x_g-\rb\right)\\
    -&\sum_g{\mbox{\boldmath $\rho$}}_g.\left( \a_{g}z_{g}-\x_g \right)\\
    -&\sum_h\Omega_h.\left(-s_h+s_{h-1}+r_h\right)\\
    -&{\mbox{\boldmath $\delta^c$}} . (t-\rb)\\
    -&{\mbox{\boldmath $\delta^d$}} . (t+\rb)\\
    -&{\mbox{\boldmath $\tau$}}.(u-\s)\\
    -&\B(\x^\frown \z^\frown t^\frown u^\frown \s)
  \end{array}
  \addtag
  \label{lagrangian}
  \]

  The marginal value identities for `room' capacity, $u$ (\ref{eq:room})
  and `door' capacity $t$ (\ref{eq:door}) are derived through the
  optimality conditions.
Particularly, when $u>0$:

  \[
\frac{\partial\mathcal{L}}{\partial u} = c^u - \sum_h\tau_{h}= 0
\implies c^u = \sum_h\tau_{h}
\]

\[
\frac{\partial\mathcal{L}}{\partial t} = c^t - \sum_h\delta_{h}= 0
\implies c^t = \sum_h\delta_{h}
  \]

(\ref{eq:room}) can be re-expressed as (\ref{eq:mvroom}) through the
  following steps. The optimality condition associated with the partial
  derivative of the Lagrangian relative to hourly storage balance produces
  the following expression, while this time explicitly incorporating the
  non-negativity constraint:
\[
s_h(\Omega_{h}-\Omega_{h+1}+\tau_h)= 0
  \]

  Combining this expression with (\ref{eq:room}), and noting that
  $\tau\geq 0$ (as the associated constraint is an inequality):
\[
\tau_h= max(0,\Omega_{h+1}-\Omega_h)\implies  c^u = \sum_h(\Omega_{h+1}-\Omega_h)^+ 
  \]

Furthermore, the implication is that $\Omega_h>\Omega_{h+1}$ can only
occur when
$s_h=0$.
The dual variable $\Omega$ thus forms a monotonically
  increasing sequence between each instant the stored energy is zero.
Summing
the differences in $\Omega$ across each sequence will lead to a
cancellation of terms, with only the first term and final term
remaining. These terms are, by monotonicity, the minimum and maximum terms
of that sequence. Indexing each such sequence, or charge/discharge cycle,
by $k$, and letting the first and final terms in the sequence $k$ be
$a(k)$ and $b(k)$ respectively, allows us to express the identity as follows:

\[
c^u=\sum_k(\max_{h\in k}\Omega_{h}-\min_{h\in k}\Omega_h)=\sum_k(\Omega_{b(k)}-\Omega_{a(k)})
\]

\subsection{What is $\Omega$?}

The optimality condition relating to the derivative of the Lagrangian with
respect to the net charge in a period allows us to investigate the
relationship between the dual variable $\Omega$ and the price of
electricity $\lambda$. The following identity shows that the two are
equivalent, when a storage unit is charging or discharging, with the
exception of a `congestion charge' on energy entering / exiting the unit.

\[
 r_h(\lambda_h-\Omega_h+\delta_{h}^c-\delta_{h}^d)=0\implies
 \Omega_h=\lambda_h+\delta_{h}^c-\delta_h^d \mid r_h\neq0
  \]

The non-negativity condition here is relevant for longer term energy
storage in particular, allowing for longer cycles between the unit
discharging to its full capacity, maintaining monotonicity of $\Omega$
through not charging or discharging in hours where the electricity price drops.

\subsection{Marginal value in general formulation}

The corresponding optimality condition in the general formulation
associated with the storage level in a
state, $s_s$, is as follows:
\[
s_s(\Omega_s-\p_s{\mbox{\boldmath$\Omega$}}+\tau_s)= 0\implies \tau_s=\p_s{\mbox{\boldmath
    $\Omega$}}-\Omega_s \mid s_s \neq 0
  \]

And by the same logic as previously it follows that:
  \[
  c^u=\sum_s(\p_s{\mbox{\boldmath $\Omega$}}-\Omega_{s})^+
  \]

  \subsection{Aggregation with guaranteed equivalence}
  \label{sec:proof-ident}

First we define aggregation for our purposes
as the partition of the set of hours into $s$ subsets. The
associated mapping between each hour $h$ and an aggregate state $s$ we term
$\Gamma(h,s)$. The weight $w_s$ of each state is the number of hours
mapped to that state $w_s=\sum_{\Gamma(h,s)}1$.
Next we define equivalence. We want the objective function value of the
aggregated model, which we denote $\bar{C}$, to equal $C$, objective
function value of the non-aggregated model. We
also want more than this typical definition, we wish the
outputs of the aggregated model, when mapped back to a disaggregated
elements, equal the outputs of the original, disaggregated, model:
$x_h=x_s$, $r_h=r_s$,
$s_h=s_s\quad\forall\Gamma(h,s)$. Our prior condition follows if this latter
condition is met.

    Recall our conditions for guaranteed equivalence:
\begin{itemize}
\item Hours that are members of a given state have equivalent temporal characteristics:
      $a_{i}=a_{j}$, $d_{i}=d_{j}$ for all hours $(i,j)$ mapped to a state $s$
\item Each state has only one incoming connection to another state, and
  only one outgoing connection: each row of the $P$ matrix contains one
  $1$ entry, with the remainder all zeros
\item For each $s$, each subsequence of hours that maps to it is of equal
  length, $q_s$
\end{itemize}

  With such a mapping, there are two steps to showing that the aggregated
  formulation will be equivalent to the non-aggregated hourly model. 
  The first step is showing that there is an optimal solution to the
  hourly model where
  the primal
  variables $x,s,r$, and dual variables $\lambda,\rho,\Omega,\delta,\tau$
in each hour $h$ mapped to the same state $s$ are equal.
This can be shown by noting that if there is an optimal solution where
hourly primal and dual values mapped to the same state are different, then
there is no barrier in terms of feasibility for them to be equal. 
This follows from the cost parameter in an objective function, $c^x$,
being constant across time, availability and demand being equal in mapped
states, and the condition of one non-zero entry of $1$ in each row of the $P$ matrix.

The second step is to apply this equality condition to the
Lagrangian (\ref{lagrangian}),
and summing the terms we develop the
following identity. We can then observe this is 
the Lagrangian of the aggregated problem (\ref{mdl:agg}).
\[
  \begin{array}{rl}
    \mathcal{\bar{L}}=&\w\sum_{g}c_g^{x}\bar{\x}_g+\c^{z}\z+c^{t}t+c^{u}u\\
    -&\bar{\lambdab}.\left(\d-\sum_g \bar{\x}_g-\bar{\rb}\right)\\
    -&\sum_g\bar{{\mbox{\boldmath $\rho$}}}_g.\left( \a_{g}z_{g}-\bar{\x}_g \right)\\
    -&\sum_s\bar{\Omega}_s.\left(-\bar{s}_s+\p_s\bar{\s}+q_s\bar{r}_s\right)\\
    -&\bar{{\mbox{\boldmath $\delta$}}} . (t-|\bar{\rb}|)\\
    -&{\mbox{\boldmath $\bar{\tau}$}}.(u-\bar{\s})\\
    -&\bar{\B}([\bar{\x}^\frown \z^\frown t^\frown u^\frown \bar{\s}])
  \end{array}
  \]
A non-obvious step in this summing of terms is to notice that
when there is a sequence of hours mapped to the same state there is a
cancelling of the $s$ terms, leaving $r$, which leads to the $q$ term
in the above identity. Also note that in the general setting we use an
overline notation to show the vectors now with dimensionality $s$ instead
of $h$.

\section{Numerical Model Description: REGEN}
\label{sec:regen}

The U.S. Regional Economy, Greenhouse Gas, and Energy (REGEN) model is used to test temporal aggregation approaches in a large-scale numerical model of long-term capacity planning with hourly operations. REGEN is fully documented in \cite{REGEN_2020}, so only a summary is provided here. The model simultaneously optimizes decisions about new generation investment, energy storage investment, hourly system dispatch, and co-optimized transmission investment and trade given assumptions about policies, technologies, and markets. The model minimizes total system costs subject to technical and economic constraints assuming perfect competition and foresight.

REGEN represents a range of energy storage technologies, including batteries, compressed air energy storage, hydrogen, and pumped hydropower. REGEN endogenously selects energy storage investment and system configurations (i.e., the ratio of energy capacity to power capacity) based on an input cost structure, where energy and power costs are specified separately. The formulation does not not include the many non-linearities associated with plant-level operations and long-run maintenance \citep{Mongrid_2019}. Battery cost and performance parameters come from \citep{Minear_2018}. Battery degradation is addressed in the model through augmentation rather than oversizing, replacement, or other strategies. REGEN includes energy storage market participation for energy arbitrage, capacity value, ancillary services (namely, operating reserves when specified), and inter-regional transmission deferral. The capacity value of energy storage and other technologies are endogenous and can vary across regions and scenarios. The formulation includes changes in investment and dispatch and endogenously considers how deployment of storage can depress its own revenues in addition to changing wholesale prices.

\section{Adjacent Clustering}\label{sec:adj}

\begin{figure}
	\centering
		\includegraphics[width=0.9\textwidth]{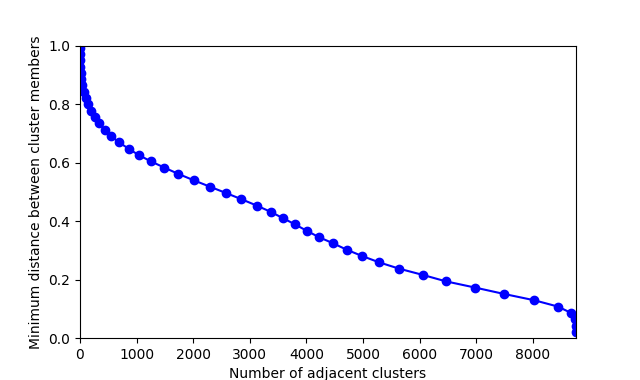}
		\caption{{Adjacent clustering
                    performance for 
                    16-region USA
                    2015 load, wind, solar dataset}}
		\label{fig:ac}
\end{figure}

  Figure \ref{fig:ac} shows the feasibility of the adjacent clustering
  approach for the REGEN USA dataset. As discussed in Section \ref{sec:general}, we require
  adjacent hours to be identical / very similar to achieve the
  conditions enabling a non-distortionary guarantee. We see in this chart that
  thousands of hours are required to achieve this theoretical guarantee, foregoing
  benefits of aggregation.

\end{document}